\documentclass[a4paper]{article}

\usepackage{amscd}

\newcommand{\marginnote}[1]{\marginpar{\tiny #1}}
\def\beq{\begin{equation}}
\def\eeq{\end{equation}}
\def\bea{\begin{eqnarray}}
\def\eea{\end{eqnarray}}

\newcommand{\Rset}{\mbox{I \hspace{-0.82em} R}}

\def\R{\Rset}

\begin{document}
\title{Introduction to the Tangent Grupoid}
\author{A. Rivero \thanks{Dept. F\'{\i}sica Te\'orica, Universidad
de Zaragoza, 50009 Zaragoza (Spain). 
Email {\tt rivero@sol.unizar.es}} }

\maketitle 

\begin{abstract}
We present some plausible definitions for the tangent grupoid of a
manifold M, as well as some of the known applications of the
structure. This is a kind of introductory note.
\end{abstract}

\vskip 1cm

Students first meet the grupoid law when they learn about
affine vectors, this is, ordered pairs of points jointly with 
an addition rule: $\vec{xy}+ \vec{yz}= \vec{xz}$. Sometimes they are able even
to notice the ambiguity of any attempt to add a "multiplication by
scalars", before being driven to quotient by an
equivalence relationship, thand then to arrive to the 
standard concept of free vector.

Perhaps the simplest mathematical object preserving the naive vector law
is the tangent grupoid. It is basically the union of two old pieces of
differential geometry: a well known one, the tangent bundle, living
in the envolvent of the manifold, and an older but less controlled
one, finite differences, which can be thought to live in the
secant envolvent of the manifold. 

 A giant step into the analysis of this structure was given by 
A. Connes, namely to study the algebra of functions over
the grupoid, using the tools of non commutative geometry \cite{connes}.

We aim 
here to present this object in a form more close to standard geometry
books (e.g.\cite{chicas,malliavin}),
so it becomes easier to focus the geometric significance of its use.

The paper can be read as a prequel to \cite[section 6]{varilly}, and we
strongly suggest the reader jump to that lecture, or some similar
one, in order to get a good grasp of the structure and its uses.               

\section{Definition}

Given a manifold $M$, and the set $A=C^1(M)$ of differentiable functions over
it, the tangent bundle $TM$ is defined as a specific subset of the
dual $A^*$: the one of evaluations of directional derivatives. To
be concrete, a vector $X$ tangent to $M$ on $x$ is associated to
the distribution (continuous linear functional) that evaluates 
at $x$ the directional derivative
along $X$:
\beq
< [x,X] | f > = \partial_X f \vert_x
\eeq

We expand this subset to a bigger one $G\subset A^*$ adding the 
distributions which give us finite differences
\beq
< [x,y,\epsilon] | f> = {f(x) - f(y) \over \epsilon} 
\eeq

%el significado de esta distribucion bien como kernel, bien
%como pre-derivada de lie, lo dejamos PARA OTRO PAPER

Even if we only had $TM$, we already confront the following problem:
While $A^*$ is not closed under products, its addition is yet an
internal operation. We would like to preserve thus this addition
in the new subset. The solution comes giving to $G$ a grupoid structure,
which restricts the pairs of elements that can be added.

We define functions range and source, from $G$ to $M \times \R$:
\beq
\matrix{
r([x,X])= (x,0) && s([x,X])=(x,0) \cr
r([x,y,\epsilon])= (y,\epsilon) && s([x,y,\epsilon])=(x,\epsilon)
}
\eeq
and the obvious immersion of $M\times \R$ in G,
\beq
e(x,\epsilon)=\left\{\matrix{ [x,0] && \epsilon = 0 \cr
                             [x,x,\epsilon] && \epsilon>0 } \right.
\eeq

With these functions, the addition in $A^*$ defines a grupoid
structure on $G$. By inducing the appropriate diferentiable
structure, $G$ can
be made a smooth grupoid, which is named {\em the tangent grupoid
to a manifold $M$}.

While the topology can be taken directly from $A^*$, for constructive
and operational reasons we can be interested on making explicit the pasting
between the two parts of the grupoid. This can be done 
\cite{connes} in an atlas
of the manifold, asking 
\beq
\lim x_n = \lim y_n = x, \lim_{\epsilon_n \to 0} {x_n-y_n\over \epsilon_n} = X.
\eeq
but other methods can be employed depending of applications.

%5 explicar (por encima) el pegado

\section{Geometric Meaning}

Before going deeper into the structure, it is a good
idea to test how it fits in our previous geometric 
visualizations. It is comfortable to check that the usual
interpretations of the tangent bundle can be enhanced to
host the full tangent grupoid. We explain this next.

\subsection{Tangent bundle and secant (pre)sheaf}

%\marginnote{grupoid OF M,  secant TO M}

When we consider $M$ immersed on some $\R^m$, there is a natural
identification of TM inside the set of straight lines in
$\R^m$, to be concrete with the ones making the tangent
envelope of $M$. This is done by seeing them as linear
applications$f:\R \to \R^m$, and then mapping each
point $[x,X]$ of $TM$ to the line tangent to $M$ in $x$ with
velocity vector $X$.

For our additional points $[x,y,\epsilon]$ there is
also an obvious subset of lines, the ones secant to $M$
through $x,y$. To keep up with the topology, we must
require that each $[x,y,\epsilon]$ is to be mapped to
the line  $f$ going through $x$ then $y$ with speed $v$ such
that ${\mbox{distance}(x,y)/|v|}=\epsilon$.
\marginnote{$\to$ We could rephrase this in terms of the natural
parameter of the straight line}

In the same manner that to each point $x \in M$ we can associate
the set of lines tangent to $x$, we can associate to every
open set $O\subset M$ the set
of lines cutting through it. While the former set is a natural
receptacle for $T_xM$, this latter one is the adequate to
map the set $S_OM\equiv \{[x,y,\epsilon] \vert x,y\in O \subset M\}$,
but the conditions to be put over the immersion for the map to
be one-to-one are more stringent in the latter case. 

Here we see one of the main weakness
of $G$. While $TM$ can be made naturally a fiber bundle
with base $M$, we have not an unique procedure to build
a fiber of elements $\{ [x,y,\epsilon] \} $ over each point
$z\in M$. The presheaf $O \to S_OM$ of elements going 
through each subset $O$ is yet unique, but it is a bit too 
general. Even if we give an arbitrary procedure to choose
a fiber, say $\pi([x,y,\epsilon])=x$, we can only get a
true vector structure in the limit $\epsilon \to 0$. 
Additionally, exotic elections for $\pi: SM\to M$ would
give problems to define the exponential map. 

%DEJAR EL RESTO DE EXP PARA OTRO PAPER

\subsection{Equivalence of curves}

Other common interpretation of $TM$ is to see each point as
a class of equivalence of curves in M, which "coincide at order $\epsilon$"
in a neighborhood of the point $x$. This
is written, given some chart $\varphi$, as:
\beq
\label{tangent.class}
[x,X]=\{ f: \R \to M \owns f(0)=x, \lim_{\epsilon \to 0} {\varphi(f(\epsilon)
) - \varphi(f(0)) \over \epsilon} = X \}
\eeq

Our generalization is obvious: each new point $[x,y,\epsilon]$ can
be seen as the class of curves passing through $x$ and $y$ for
determinate values of its parameter. For example, if we want
to be consistent\footnote{The ambiguity in the definition of derivative
must translate to an ambiguity in the value of the parameter of the
curve which goes through $x$.} with 
the previous formula for $TM$, we can postulate:
\beq
\label{secant.class}
[x,y,\epsilon]=\{f: \R \to M \vert f(0)=x, f(\epsilon)=y \}
\eeq

Both spaces are pasted using simple set-theory techniques: a sequence
$\{ [x_n,y_n,\epsilon_n] \}$ is said to converge to $[x,X]$ 
if and only if $\{\epsilon_n\}$ goes to zero and
\beq
\emptyset \not= \bigcap_n [x_n,y_n,\epsilon_n] \subset [x,X]
\eeq

Alternatively, we can even avoid to define $TM$ and postulate it
directly as the class of non empty intersections of
such sequences.

\subsection{Leibnitz rule}

A vector of $T_xM$ can be defined\footnote{In finite dimensional
 manifolds \cite{chicas}.} also as
an element $v_x \in A^*$
fulfilling Leibnitz rule, 
\beq
v_{x}(fg)= f(x) v_{x}(g) + g(x) v_{x}(f)
\eeq

The definition can be expanded to
the rest of the grupoid by permitting "braided" Leibnitz rules:
\beq
\label{leibnitz.braided}
v_{xy}(fg)= f(x) v_{xy}(g) + g(y) v_{xy}(f)
\eeq

In this sense, distributions $\left<[x,y,\epsilon]\right|$ are 
considered as "deformed" derivations. Sometimes formula 
(\ref{leibnitz.braided}) is written using coordinates more or
less implicitly; e.g. with a coordinate function $\phi$:
\def\smphi{{}_{\tiny \phi}}
\beq
v_{xy}(f)={f(x)-f(y)\over \smphi x-\smphi y} v_{xy}(\phi)
\eeq 

\subsection{TM as boundary of the grupoid}

The set-theory building of $TM$ as limiting set of $SM$
was already apparent in section 2.1, where we can observe 
than the set of tangent lines lies in the border
of the set of secant ones. 

The goal of pasting the two parts of $G$ in a form compatible
with the usual definition of derivatives brings this result
as byproduct: the tangent bundle is the boundary of the
secant grupoid.

%%% NO HABLAR DEL DOUBLING EN ESTE PAPER; DEJARLO PARA OTRO!!!!

We take usually $\epsilon\in [0,1)$ when we want to see G as a manifold 
with boundary. In this case we will conventionally call $G^{(0)}\equiv TM$ to
the tangent part, $\epsilon=0$, and $G^{(1)}\equiv SM$ to the secant part,
$0<\epsilon<1$. 

\section{A Dilatation Structure}

The pasting of $SM$ to $TM$ uses some charting of $M$ in order to
control the rate of convergence of $x_n,y_n$. This is, to say the
least, antieconomical, as only the scaling structure provided
by the chart is used. In addition, the atlas itself is unavailable if 
we choose to define TM directly from sequences in $SM$.

Alternatively, we can control the convergence by defining a
dilatation structure over SM
\beq
\tau_\lambda:SM \to  SM 
\eeq
\beq
[x,y,\epsilon] \to [x',y',\lambda \epsilon] 
\eeq
%\marginnote{FALTA ALGUNA CONDICION PARA CONTROLAR LA LINEALIDAD}
with $\tau_{\lambda \lambda'}=\tau_\lambda \circ \tau_{\lambda'}$
(For a general discussion on this structure, see \cite[ch.12]{W-K}).

We will say that a sequence $\{[x_n,y_n,\epsilon_n]\}$ going to
$[x,x,0]$ has a limit on $TM$ ("goes to a point of $TM$") if the scaled
 sequence $\{\tau_{\varepsilon_0
/ \epsilon_n}[x_n,y_n,\epsilon_n]\}$ has a limit in $SM$.
\marginnote{($SM= M\times M \times \R^+$)}

%\marginnote{Debemos preservar el producto del grupoide al escalar? O
%solamente en el limite renormalizado?}

In some sense, this technique is reminiscent of the tensorial definition:
"a vector is a geometrical object which transforms as a vector", but here
transformations are hidden under the carpet of the dual $A^*$.

If we need explicitly to know the limit vector on the tangent bundle,
we must do a bit more of work, as remarked in section one, or, alternatively,
specify the recipe compatible with the dilatation flow. This is not surprising:
we have used the flow as a substitute for the explicit chart in
(\ref{tangent.class}), but then we have lost a method to name
explicitly $X$. We must then specify a recipe compatible with the
"hidden" coordinate system. 

For example, take the flow $\tau_{\lambda}$ to be such that
\beq
\label{example.bohr}
\tau_\lambda [x, \exp_x X, \epsilon] = [x, \exp_x \lambda X, \lambda \epsilon]
\eeq
Then each point $[x,X]$ must be canonically associated to the
limit point $[x, \exp_x X, \varepsilon_0]$.

To put other common example, if we choose a "mid-point rule" for
scaling:
\beq
\label{example.weyl}
\tau_\lambda [\exp_z -X, \exp_z X, \epsilon] =
[\exp_z -\lambda X, \exp_z \lambda X, \lambda \epsilon]  
\eeq
the limit point $"[x,Y]"$ will be $[\exp_x -Y, \exp_x +Y, \varepsilon_0]$
 
%aqui nos debe permitir definir TM como este limite, sin recurrir
% a misteriosas derivaciones delta' en A*

%when works, when no? Ok for Lie derivative, but not for covariant?

We have not a general proof of existence for this process 
of limit. It could be suspected that, while working well for
directional derivatives, it could fail for more exotic
constructions, such as covariant derivatives depending 
% o simplemente para hacer conexiones (un campo gauge es una
% conexion en un fibrado principal, una derivada covariante
%es su aplicacion en un fibrado vectorial asociado)
on symmetry groups (from a physicist point of view, we would like
geometry to be able to detect physical properties of the 
connection, specifically renormalizability).

%explicacion de la hipotesis de trabajo, no se imprime en el texto:
% el criterio general es que el hecho de que una conexion sea o no
% una TCC renormalible deberia aparecer en algun sition en la
% geometria.
% la hipotesises es que el sistema de escalado depende de tener puntos
% fijos. Para la derivada normal el punto fijo correspondiente nos da
% la definicion usual de derivada, digamos que es una q-derivada
% con q->1. Para conexiones mas complicadas puede que el limite
% continuo tenga que apoyarse en puntos fijos mas complicados, y que
% este se refleje en derivadas donde q debe tender a otras raices
% de la unidad. A continuacion esperariamos que el limite pudiera
% ser descrito de forma aproximada como la suma de varias cosas
% con q->1. En el caso de interes, para gravedad, la conexion
% asociada necesitaria una raiz por determinar (quizas cubica)
% y podria aproximarse como una suma de los casos q->1 de la
% misma conexion y de la conexion asociada al modelo standard.
%
% si contamos explicitamente el numero de terminos diferentes de la
% accion de gravedad, y su agrupacion en tipos segun el numero
% de derivadas que implican, vienen a se el mismo numero de libertades
% que para el algebra que codifica el modelo standard. Esto podria
% ser una pista. 

\section{The Tangent Grupoid of a Lie Group}

When the manifold $M$ has a Lie group structure, we can say
a lot more about the tangent bundle: we know how to get each fiber
$T_xM$ from an action of the group over the fiber at the identity,
$T_eM$. And we know that this fiber, the tangent space at
the identity, can be show to be equivalent to the Lie algebra {\cal g}
of the group (and/or to the left invariant fields over M), and we can
build $A^*$ from the enveloping algebra of {\cal g}.

Now, if we are able to get such quantity of information
from the $\epsilon=0$ part of our grupoid, it is natural to ask
ourselves how this information is present in any part of 
it $\epsilon=\epsilon_0
>0$. A partial answer in the language of
quantum groups is provided 
by \cite{majid}, which builds the distributions associated to the
tangent grupoid as a particular first-order example of more general
deformed (braided) differential calculus over M. There, a explicit selection
of fiber over each $x\in M$ is needed, entering the yet developing 
world of deformed q-vector bundles.

The theory of differential calculus for q-groups is currently
well established. Basically, it deforms both the algebra $A$ of functions over
$M$ and the dual algebra of distributions $A^*$ (to be precise, the
enveloping algebra, as we have said), but keeps both in duality by using 
(an equivalent of) the usual relation
\beq
\label{duality}
<X, f> =    Xf |_0, \,\, f\in A, X \in g \subset A^* 
\eeq
Some additional structure is required, mainly as properties of
the coproducts, see \cite{Comm}. In the clasical, undeformed, case,
these are $\Delta(X)= X\otimes 1 + 1 \otimes X$ and
$\Delta(a(x))=a(x . y)$, for $X\in A^*$ and $a\in A$ respectively.
\marginnote{ currently, Sweedler notation
is used, $X_{(1)}\otimes X_{(2)}$ and $a_{(1)}\otimes a_{(2)}$
respectively}

The tangent space to the identity, L, is characterized technically through
the action of the quantum
double $D(A^*)$ in $\mbox{ker } e \subset A^*$, where 
$e$ is the counit of $A^*$. The usual product of functions
and vector fields, $a(x)X$, generalizes to the action
\beq
a \triangleright y = <a, y_{(1)}> y_{(2)} -1 <a,y>
\eeq
while the Lie derivation of vector fields retain its aspect as
adjoint action $Ad_x$ 
\beq
x \triangleright y = x_{(1)} y S x_{(2)}
\eeq
\marginnote{$X \triangleright Y= XY - YX $} 

In a dual manner, over $A$ we can locate a space $L^*=\mbox{ker} e' /M$, being 
M some $D(A)$-submodule and $e'$ the counit of A. Such space will
have a natural action of $A$ by multiplication from the left, and an
acting of $A^*$ as coadjoint action, $Ad^*$, both following from the
construction of the double $D(A)$. These actions are used to define the
space of left-invariant 1 forms over the group, explicitely 
$\Omega^1=L^* \otimes A$. 

%\beq
%a.(v\otimes b) = a_1 \triangleright v \otimes a_2 b
%\eeq
%(where $\triangleright$ is (left) multiplication and the projection over 
%$L^*$) and
%\beq
%(v \otimes b) .a = v \otimes b a
%\eeq
%\beq
%da=(\pi \otimes id) (1 \otimes a - a_1 \otimes a_2)
%\eeq
%(where $\pi$ is the canonical projection over $L^*$ we have used.)
%
%\marginnote{podriamos dar el ejemplo de la projeccion en nuestro caso
%concreto}
%...
%

Details of actions on $\Omega^1$ are show in \cite{majid} jointly
with 
explicit examples for $\R^n$. In particular, when
$n=1$ any function $c(p)$ defines an invariant q-tangent space
\beq
L=\mbox{span } \{p_n=c^{(n)}(p)-c^{(n)}(0)\}
\eeq
with corresponding (braided) derivation and Leibnitz rule, which
we omit for sake of brevity.

We obtain $L^*$ taking the quotient of $\mbox{ker} \epsilon'$, this is,
functions $f \owns f(0)=0$, by its subspace
\beq
\{ f \owns  \partial_{p_n}f |_0 =0 \}
\eeq
thus keeping with the usual method given by
(\ref{duality}).

For $c_\lambda(p)=\lambda^{-2} e^{\lambda p}$, L is one-dimensional, and
its partial derivative comes out to be 
\beq
\label{partial}
\partial_{p_1} f = {f(x+\lambda) - f(x) \over \lambda}
\eeq
while for $c_{\lambda \to 0}$ we get the usual commutative, $q\to 1$,
calculus, which is also our limiting process $\epsilon\to 0$. For other
derivatives it would be interesting to investigate limits to other
roots of the unit.

\section{Applications}

Connes approach \cite{connes} studies the tangent grupoid through its
Gelfand-Naimark dual, i.e. the set of continuous (alt. smooth)
functions over it, which form an algebra, its product being
the convolution product associated to the grupoid law.
\beq
(a*b) (\gamma) = \sum_{\matrix{\alpha,\beta \in G \cr \alpha\beta=\gamma}}
		     a(\alpha) b(\beta)
\eeq

For fixed $\epsilon$, greater than zero, this is product of
kernels of operators in $L^2(M)$ (see \cite{varilly}). Thus
we have is a very nice host for a K theory.

In fact the main use of the tangent grupoid has been to
work its K theory as receptacle for Bott periodicity, then
to formulate proofs of index theorems.

Also, it has been observed \cite{houches} that the set of continuous functions
over the grupoid defines a quantization procedure. This is
because a observable function over $T^*M$ can be described
with its Fourier transform in $TM$, the pointwise product passing
then to convolution product. If we assume that the function in
TM is the $\epsilon = 0$ part of a continuous function over the
full grupoid, it can be seen that any fixed $\epsilon=\epsilon_0$ part
of it is described as the kernel of a quantized operator associated
to the original function. The ambiguity of this procedure
is studied in \cite{ours}

%tambien remarcar la formula rara de Majid aqui al final, y no antes

It could be worth to note that the quantization phenomena
appears also in the example of \cite{majid}, where formula (19)
can be seen as implementing the indeterminacy principle needed
in any quantization context (If we prefer to see $A$ as a set of
fields on 0+1 dimensions, then it is most appropriate to try interpret (19)
as the quantum mechanics evolution rule).

\section*{Acknowledges}

Author must acknowledge Jose Gracia-Bond\'{\i}a insistence to induce
our team to work in the grupoid tangent, and Jes\'us Clemente, J.F.
Cari\~nena 
and Eduardo Follana for discussions and critical reading. The appendix
specially has benefited from comments from Eduardo, and the author must
apologize for not implementing all the suggested improvements. 

%a los lectores por su paciencia 
% dftuz: alojamiento

%\newpage

\section*{Appendix: Toy Renormalization Group Triangle}

Lets review how we use a dilatation structure on our set of distributions
$\{\left<[x,y,\epsilon]\right|\}$. Ideally, this structure is a set of
transformations $\{\tau_\lambda\}$ which are obtained by duality 
from other set  $\{\tau_\lambda^*\}$ of "rescaling" 
operators on $A=C^\infty(M)$. 
 \beq
\label{rg.scale}
<\tau_\lambda v | \tau^*_\lambda a> = < v | a> 
\eeq
Such pairing can be built explicitly in simple cases, by example for 
$\R$ we can put $\tau_\lambda$ as in (\ref{example.bohr}) and $\tau^*\lambda$ 
such that for some family of charts $\phi$, any function $a \in A$ transforms 
as in 
a linear change of variable,
\beq
\tau^*_\lambda a (x) = \lambda^1 a(\phi^{-1} \lambda^{-1} \phi x).
\eeq
But regretfully the pairing is not one on one, so in practice
(\ref{example.weyl}, \ref{example.bohr}, etc.)  we start
directly from a concrete $\{\tau_\lambda\}$ set. 

Now, dilatations $\tau_\lambda$ draw a flow on the
space $(M\times M \times \R)$ that parametrizes $SM\subset A^*$, and
this 
can be interpreted as a RG flow. Lets see how.

We put some mild restrictions in transformations $\tau_\lambda$. They
must form a multiplicative semigroup, 
$\tau_{\lambda\mu}=\tau_\lambda \tau_\mu$. They must preserve the subspace
$M\times M\times \R$ of distributions. In order to simplify operations, they
shall act multiplicatively in $\R$, 
$\tau_\lambda [x,y,\epsilon] = [x',y',\lambda\epsilon]$

The dilatation condition takes apart any pair
of points $x,y$, and additionally moves $\epsilon$. Then the only
fixed points of the flow are $\{[x,x,0]\}$.
The standard approach ask us to linearize the transformation
near the fixed point, and then study its flow.
 
The fixed point lies in the critical surface 
$\{ \left<[x,y,0]\right|, x,y \in M \}$, where distributions become
indefinite. Any
series going to a point $[z,z,0]$ in the critical surface 
draws a nonrenormalized line of "bare" distributions
\beq
\left\{[x_n,y_n,\epsilon_n]\right\}, \epsilon_n\to 0,
 y_n \to x_n \to z, n\to\infty
\eeq
To fix it, we choose a "physical" scale $\varepsilon_0$ and use
the renormalization flow transformations $\tau_{\varepsilon_0/\epsilon_n}$
to rescale the points, getting a renormalized series
\beq
\left\{ \tau_{\varepsilon_0\over\epsilon_n} [x_n,y_n,\epsilon_n]\right\}
\eeq
which lives in the surface $[x,y,\varepsilon_0], x,y\in M$. If the
bare series had a limit $[x,X]$ in the associated chart $\phi$, then
the renormalized series has a limit, which we can name
$[x,X]_{\varepsilon_0}$.

The mathematical image is straightforward: both limits are
related "until energy $\varepsilon_0$", which is the scale of
validity of the renormalized one;
\beq
\langle[x,X]| a\rangle = \left<[x,X]_{\varepsilon_0}|a\right> + o(\varepsilon_0)
\eeq
Renormalized limits at different scales are joined by RG transformations,
correspoding to a relevant direction coming out from the fixed point.

Other usual folklore of Renormalization Group theory can equally be
translated to this context: So, formula (\ref{rg.scale}) is the
traditional assesment which claims physics to be independent of the
scale; the arbitrarity choosing an specific R.G. transformation
translates to the freedom to choose the limiting process that defines
the derivative of a function, while the global character of the
obtained tangent vector comes from the global properties of
the fixed point and the flow near it.

The main drawback of the example is the character of the fixed point:
it has not got any irrelevant direction, then any "bare" line
must go directly to it in order to define a good limit.
 This effect can be attributed to the
simple form of the derivative being defined in the process. If we want to look
for examples owning irrelevant directions in the critical surface, we would 
try more complicated geometrical notions, such as covariant derivatives.
%quizas las irrelevantes sean las elecciones de horizontal lifts, ver
% 10.49 en el nakahara. O quizas dependan directemente del vector
% que nos interese derivar covariantemente
Also, we could try Jackson q-derivatives, for $q$ going to a complex
root of unit.

%de hecho, ambas cosas estaran relacionadas. Algunas derivadas covariantes
% estaran definidas solo para ptos fijos peculiares, que corresponderan
% a f(qx)-f(x) / (q^N-1)x, o a su representacion finito dimensional o
% como sea que va eso.
% otro punto de misterio es si tales ptos fijos deberian verse como
% "antiferromagneticos", ya que son limites continuos muy peculiares.

\end{document}